\documentclass[]{article}
\usepackage{graphicx}
\usepackage{graphics}
\usepackage{amssymb,amsmath,amsthm,epsfig}
\newtheorem{theorem}{Theorem}
\newtheorem{lemma}{Lemma}

\newtheorem{proposition}{Proposition}

\newtheorem{remark}{Remark}
\newcounter{rea}
\newcommand{\wJ}{\widetilde P^{(\alpha,\beta)}}
\newcommand{\J}{P^{(\alpha,\beta)}}
\newcounter{rek}
\newcommand{\ps}{\psi_{n,c}^{(\alpha)}}
\newcommand{\pn}{\widetilde P_n^{(\alpha,\alpha)}}
\newcommand{\qn}{\widetilde Q_n^{(\alpha,\alpha)}}
\begin{document}

\begin{center}
{\large {\bf Weighted finite Fourier transform operator: Uniform approximations of the eigenfunctions, eigenvalues  decay  and behaviour.}}\\
\vskip 1cm Abderrazek Karoui$^a$  {\footnote{
Cooresponding author: Abderrazek Karoui, Email: abderrazek.karoui@fsb.rnu.tn
This work was supported by the DGRST research Grant UR13ES47 and the CMCU Research project 15G 1504.}} and Ahmed Souabni$^a$
\end{center}
\vskip 0.5cm {\small

\noindent $^a$ University of Carthage,
Department of Mathematics, Faculty of Sciences of Bizerte, Tunisia.
}\\

\noindent{\bf Abstract}--- In this paper, we first give two uniform asymptotic approximations of the 
eigenfunctions of the  weighted finite Fourier transform operator, defined by ${\displaystyle \mathcal F_c^{(\alpha)} f(x)=\int_{-1}^1 e^{icxy}  f(y)\,(1-y^2)^{\alpha}\, dy,\,}$ where $ c >0, \alpha > -1$ are two fixed real numbers.
The first uniform approximation is given in terms of a Bessel function, whereas the second one is given in terms
of a normalized Jacobi polynomial.
These eigenfunctions are called generalized prolate spheroidal wave functions (GPSWFs). By using the uniform asymptotic approximations of the GPSWFs,  we prove the  super-exponential decay 
rate of the eigenvalues  of the operator $\mathcal F_c^{(\alpha)}$ in the case where $0<\alpha < 3/2.$
Finally, by computing   the trace and an estimate of the norm of the operator ${\displaystyle \mathcal Q_c^{\alpha}=\frac{c}{2\pi}
\mathcal F_c^{{\alpha}^*}  \mathcal F_c^{\alpha},}$ we give a lower and an upper bound for the counting number of the eigenvalues of $Q_c^{\alpha},$ when $c>>1.$\\



\section{Introduction} In the early 1960's, D. Slepian and his co-authors H. Landau and H. Pollack, have greatly contributed in developing the theory of prolate spheroidal wave functions (PSWFs), see their  pioneer work \cite{Landau, Slepian1, Slepian2, Slepian3}. For $c>0,$ a positive real number, called bandwidth, the 
PSWFs, denoted by $(\psi_{n,c})_{n\geq 0}$ are the eigenfunctions of the finite Fourier transform operator ${\displaystyle \mathcal F_c }$, as well as the Sinc kernel convolution operator 
${\displaystyle \mathcal Q_c,}$ defined on $L^2([-1,1])$ by  
${\displaystyle \mathcal F_c f(x)=\int_{-1}^1 e^{icxy}  f(y)\, dy,\quad  \mathcal Q_c (f)(x)=\frac{2c}{\pi} \mathcal F_c^{*} \circ \mathcal F_c.}$ Perhaps the starting point of the theory of PSWFs is the Slepian's result concerning the commutativity property of the integral operators  ${\displaystyle \mathcal F_c }$ and ${\displaystyle \mathcal Q_c }$ with the following perturbed Legendre differential operator
$$\mathcal L_c y(x)= -(1-x^2) y"(x) + 2 x y'(x) + c^2 x^2 y(x).$$
Since $L_c \mathcal F_c = \mathcal F_c L_c,$ then the PSWFs are also, the bounded eigenfunctions over  $I=[-1,1]$ of the Sturm-Liouville operator $\mathcal L_c.$ Many desirable properties, computational schemes, asymptotic results and expansions of the PSWFs are consequences of the previous commutativity property. This important property  has allowed the use and the application of the rich literature of the theory 
of Sturm-Liouville operators in the context of the PSWFs. 

We should mention that the PSWFs have found applications in various  area such as 
applied mathematics, mathematical physics, random matrices, signal processing, etc., see \cite{Hogan}
for a comprehensive review of the theory and some applications of the PSWFs. Note that most of the PSWFs applications, rely of the decay rate and the behaviour of the eigenvalues of the integral operator $\mathcal F_c$ or of $\mathcal Q_c,$
as well as the bounds and the local estimates of the PSWFs. 

Recently, there is an interest in the spectral analysis of a more general compact integral operator, 
the weighted finite Fourier transform operator $\mathcal F_c^{(\alpha)},$ defined  by
\begin{equation}\label{Eq1.1}
\mathcal F_c^{(\alpha)} f(x)=\int_{-1}^1 e^{icxy}  f(y)\,\omega_{\alpha}(y)\, dy,\quad \omega_{\alpha}=(1-y^2)^{\alpha},\quad \alpha >-1.
\end{equation}
It is well know, see \cite{Karoui-Souabni1, Wang2} that the operator ${\displaystyle \mathcal Q_c^{\alpha}=\frac{c}{2\pi}
\mathcal F_c^{{\alpha}^*}  \mathcal F_c^{\alpha}}$ is  defined on $L^2{(I, \omega_{\alpha})}$ by
\begin{equation}
\mathcal Q_c^{\alpha} g (x) = \int_{-1}^1 \frac{c}{2 \pi}\mathcal K_{\alpha}(c(x-y)) g(y) \omega_{\alpha}(y) \, dy,
\quad \mathcal K_{\alpha}(x)=\sqrt{\pi} 2^{\alpha+1/2}\Gamma(\alpha+1) \frac{J_{\alpha+1/2}(x)}{x^{\alpha+1/2}}.
\end{equation}
The eigenvalues $\mu_n^{(\alpha)}(c)$ and $\lambda_n^{(\alpha)}(c)$ of $\mathcal F_c^{\alpha}$ and 
$\mathcal Q_c^{\alpha}$ are related to each others by the identity
${\displaystyle \lambda_n^{(\alpha)}(c)= \frac{c}{2\pi } |\mu_n^{(\alpha)}(c)|^2}.$ 
Moreover, both operators commute with the following Jacobi-type Sturm-Liouville operator $\mathcal L_c^{(\alpha)},$ defined by
$$
\mathcal L_c^{(\alpha)} (f)(x)= - \frac{d}{dx}\left[ \omega_{\alpha}(x) (1-x^2) f'(x)\right] +c^2 x^2 \omega_{\alpha}(x) f(x).$$
The infinite countable set of the eigenfunctions of $\mathcal F_c^{\alpha}, \mathcal Q_c^{\alpha}$ and $\mathcal L_c^{(\alpha)}$ will be denoted by $(\ps)_{n\geq 0}.$ They  are called generalized prolate spheroidal wave functions 
(GPSWFs). Some properties as well as a first set of local estimates and bounds of the GPSWFs have been given in  
\cite{Karoui-Souabni1}. 
It can be easily checked that $\chi_{n,\alpha},$ the  $n-$th eigenvalue of the differential operator $\mathcal L_c^{(\alpha)}$ satisfies the inequalities, see \cite{Karoui-Souabni1}
$$n(n+2\alpha+1)\leq \chi_{n,\alpha}\leq n(n+2\alpha+1)+c^2,\quad n\geq 0.$$
Also, it has been shown in \cite{Bonami-Karoui1} that in the special case where $\alpha=0,$ the eigenvalues 
$\lambda_n(c)=\lambda_n^{(0)}(c)$ decay asymptotically faster than  $e^{-2n\log\left(\frac{a n}{c}\right)}$ for any positive real number $0<a < \frac{4}{e}.$  In \cite{Wang2}, for a more general value of 
$\alpha >-1,$ the authors have checked
that for the sequence of the eigenvalues $\lambda_n^{(\alpha)}(c)$ have an asymptotic decay rate similar to the sequence
$e^{-(2n+1) \log\left(\frac{4n+4\alpha+2}{e c}\right)}.$ Nonetheless, this result is obtained by using some  heuristic results concerning the behaviour and the decay of the coefficients of the Gegenbauer's series expansion of $\ps.$ 
In this work, we give a proof of the previous super-exponential decay rate of the $\lambda_n^{(\alpha)}(c)$ in case where $0<\alpha < 3/2.$ This proof is based on two uniform asymptotic approximations of the $\ps.$ The first one is given in terms of the Bessel function $J_{\alpha}(\cdot)$ and the second one is given in terms of the normalized 
Jacobi polynomial $\pn.$  Note that the $\ps$ and $\pn$ are normalized by the following rules
\begin{equation}\label{normalisations}
\int_{-1}^1 (\pn(x))^2 \omega_{\alpha}(x)\, dx = 1,\quad \int_{-1}^1 (\ps(x))^2 \omega_{\alpha}(x)\, dx = 1.
\end{equation}
Under the above normalisation of $\ps,$  we show that 
for any positive integer $n$ with  ${\displaystyle q=\frac{c^2}{\chi_{n,\alpha}} <1,}$ we have
$$\ps \approx   \sqrt{\frac{\pi}{2\mathbf K (\sqrt{q})}} \frac{(\chi_{n,\alpha})^{1/4}\sqrt{S(x)}J_{\alpha}(\sqrt{\chi_{n,\alpha}} S(x))}
{(1-x^2)^{1/4+\alpha/2}(1-q x^2)^{1/4}}, \quad x\in [0,1].$$
Here, 
${\displaystyle S(x)=\int_x^1 \sqrt{\frac{1-qt^2}{1-t^2}}\, dt,\,\,\, \mathbf K(r)=\int_0^1 \frac{1}{\sqrt{(1-t^2)(1-r^2 t^2)}},\, dt \,\,\, 0\leq r<1.}$  Also, by using some properties and estimates of the Jacobi polynomials and Jacobi functions of the second kind, we prove that if $0<q \leq q_0<1,$ then  for sufficiently large values of $n,$ we have the following uniform approximation of the $\ps$ in terms of the Gegenbauer's polynomial,
$$\left|\ps(x)-A_n \pn (x)\right|\leq C_{\alpha}(q_0) \frac{  c^2}{n+2\alpha+1},\quad \forall\,\, x\in [-1,1],
$$
where  $A_n$ is a normalization constant, satisfying ${\displaystyle |1- A_n | \leq C_{\alpha}(q_0) \frac{c^2}{2n+2\alpha+1}}$ and $C_{\alpha}$ is a constant depending only on $\alpha.$ 
Also, we show that  if  $ 0<\delta<1$ and if  $M_c(\delta)$ is  the number of eigenvalues of $Q_c^{(\alpha)},$  $ \alpha>0$, which are not smaller than $ \delta,$ then
 \begin{equation}
 \frac{\gamma_{\alpha}-\delta}{1-\delta}\frac{c}{2\pi} (2^{2\alpha+1}B(\alpha+1,\alpha+1))^2 +o(c)\leq  M_c(\delta) \leq \frac{1}{\delta}\Bigg[ \frac{c}{2\pi} \Big[ 2^{2\alpha+1} B(\alpha+1,\alpha+1) \Big]^2 \Bigg].
 \end{equation}
 Here, $ \gamma_{\alpha}=2^{4\alpha}\Big( \frac{B(2\alpha+1,2\alpha+1)}{B(\alpha+1,\alpha+1)} \Big)$ and $B(\cdot,\cdot)$ is the beta function.\\

This work is organized as follows. In section 2, we study a uniform asymptotic approximation of the 
$\ps$ in terms of the Bessel function $J_{\alpha}.$ This approximation result is based on the use of 
the  WKB and Olver's methods, together with some properties of Bessel functions. In section 3, we first list  some properties
and estimates of the Jacobi polynomials and Jacobi functions of the second kind. Then by using these results, we prove
the asymptotic uniform approximation of the $\ps$ in terms of the Jacobi polynomials $\pn.$ 
In section 4, we first use the result of the previous two sections and prove the super-exponential decay rate of the 
$\lambda_n^{(\alpha)}(c),$ for $0<\alpha< 3/2.$ Then, by using the trace and an estimate of the norm of the integral
 operator $\mathcal Q_c^{\alpha},$ we give a  lower and an upper  bound for the counting number of  the
 eigenvalues $\lambda_n^{(\alpha)}(c).$

\section{Uniform approximation of the eigenfunction in terms of Bessel functions.}
Let $ w_{\alpha}(x)=(1-x^2)^{\alpha}$ and recall that 
 the GPSWFs are also the bounded eigenfunctions on $I=[-1,1]$ of the following differential equation,
\begin{equation}\label{eq1_diff}
(1-x^2)\psi''(x)-2(\alpha+1)x\psi'(x)+(\chi_{n,\alpha}-c^2 x^2)\psi(x)=0.
\end{equation}
Here, $\chi_{n,\alpha}$ is the $(n+1)$th eigenvalue of the following Sturm-Liouville differential operator
\begin{equation}\label{Sturm_form}
\mathcal L_c^{(\alpha)} (f)(x)= - \frac{d}{dx}\left[ \omega_{\alpha}(x) (1-x^2) f'(x)\right] +c^2 x^2 \omega_{\alpha}(x) f(x)
\end{equation}
Recall that the eigenvalue $\chi_{n,\alpha}$ satisfies the following classical inequalities,
\begin{equation}\label{boundschi}
n (n+2\alpha+1) \leq \chi_{n,\alpha} \leq n (n+2\alpha+1) +c^2,\quad \forall n\geq 0.
\end{equation}
The previous differential equation is rewritten as 
\begin{equation}\label{eq2_diff}
-\mathcal L_c^{(\alpha)} \psi(x)+w_{\alpha}(x)\chi_{n,\alpha}\psi(x)= ( w_{\alpha}(x)\psi'(x)(1-x^2))'+ w_{\alpha}(x)(\chi_{n,\alpha}-c^2 x^2)\psi(x)=0,\quad x\in [-1,1].
\end{equation}
We use the Liouville transformation to transform this later equation into a Liouville normal form. 
More precisely, for a positive integer $n$ with ${\displaystyle q=\frac{c^2}{\chi_{n,\alpha}} <1,}$ we consider 
the incomplete elliptic integral 
\begin{equation}\label{Liouvile_transform1}
S(x)=\int_x^1 \sqrt{\frac{1-qt^2}{1-t^2}}\, dt.
\end{equation}
It has been shown in \cite{Bonami-Karoui2},  that for  $0< q <1,$ we have
\begin{equation}\label{equS}
(1- \frac q2) \sqrt{(1-x^2)(1-q x^2)} \leq S(x) \leq \frac{5-q}3\sqrt{(1-x^2)(1-q x^2)}.
 \end{equation}
Then, we  write $\ps$ into the form
\begin{equation}\label{Liouvile_transform2}
 \psi(x)=\phi_{\alpha}(x)V(S(x)),\qquad \phi_{\alpha}(x)=(1-x^2)^{(-1-2\alpha)/4} (1-qx^2)^{-1/4}.
\end{equation}
By combining (\ref{eq2_diff}), (\ref{Liouvile_transform2}) and using straightforward computations,
it can be easily checked that $V(\cdot)$ satisfies the following second order differential equation 
\begin{equation}
V''(s)+\left(\chi_{n,\alpha}+\theta_{\alpha}(s)\right)V(s)=0,\quad s\in [0, S(0)]
\end{equation}
with $$ \theta_{\alpha}(S(x))=(w_{\alpha}(x)(1-x^2)\phi_{\alpha}'(x))'\frac{1}{\phi_{\alpha}(x)w_{\alpha}(x)(1-qx^2)}.$$
If $ Q_{\alpha}(x)=w_{\alpha}(x)^2(1-x^2)(1-qx^2),$ then we have 
${\displaystyle  \frac{\phi_{\alpha}'(x)}{\phi_{\alpha}(x)}=-1/4 \frac{Q_{\alpha}'(x)}{Q_{\alpha}(x)}. }$
It follows that $\theta_{\alpha}(S(x))$ can be written as
\begin{equation}
\theta_{\alpha}(S(x))=\frac{1}{16(1-qx^2)}\Big[\left(\frac{Q_{\alpha}'(x)}{Q_{\alpha}(x)}\right)^2(1-x^2)-4\frac{d}{dx}\Big((1-x^2)\frac{Q_{\alpha}'(x)}{Q_{\alpha}(x)}\Big)-4(1-x^2)\frac{Q_{\alpha}'(x)}{Q_{\alpha}(x)}
\frac{w_{\alpha}'(x)}{w_{\alpha}(x)}\Big].
\end{equation}
Since $ Q_{\alpha}(x)=w_{\alpha}^2(x)Q_0(x),$ then we have 
${\displaystyle  \frac{Q_{\alpha}'(x)}{Q_{\alpha}(x)}=2\frac{w_{\alpha}'(x)}{w_{\alpha}(x)}+\frac{Q_0'(x)}{Q_0(x)}} $ and ${\displaystyle \frac{w_{\alpha}'(x)}{w_{\alpha}(x)}=-\frac{2\alpha x}{1-x^2}.} $
Hence, we have 
\begin{eqnarray}
\theta_{\alpha}(S(x)) &=&\theta_0(S(x))+\frac{-1}{4(1-qx^2)}\Big[ \big(\frac{w_{\alpha}'(x)}{w_{\alpha}(x)}\big)^2 (1-x^2)+2\frac{d}{dx}[(1-x^2)\frac{w_{\alpha}'(x)}{w_{\alpha}(x)}]  \Big] \\
&=&\theta_0(S(x))+ \frac{1}{(1-x^2)(1-qx^2)}\Big( -\alpha^2 x^2 + \alpha (1-x^2) \Big) \\
&=&\theta_0(S(x))+ \frac{\alpha(1+\alpha)}{(1-qx^2)}-\frac{\alpha^2}{(1-q x^2)(1-x^2)}.
\end{eqnarray}
It has been shown in \cite{Bonami-Karoui2} that 
$$\left|\theta_0(x)-\frac{1}{4 S^2(x)}\right|\leq \frac{3+2q}{4(1-q x^2)^2},\quad \left|\frac{1}{S^2(x)}-\frac{1}{(1-q x^2)(1-x^2)}\right|\leq \frac{3}{(1-q x^2)^2},\quad  x\in (-1,1).$$
Consequently, if $G_{\alpha}(\cdot)$ is the function given by 
\begin{equation}
\label{fctG}
G_{\alpha}(x)=\frac{1/4 -\alpha^2}{S^2(x)}-\theta_{\alpha}(x),\quad x\in (-1,1),
\end{equation}
then, we have
\begin{equation}
\label{bound1fctG}
|G_{\alpha}(x)|\leq \frac{3+2q+12 \alpha^2}{4(1-q x^2)^2}+\frac{\alpha(\alpha+1)}{1-q x^2},\quad x\in (-1,1),
\end{equation}
As it is done in \cite{Bonami-Karoui2}, by using the substitution $t=S(x),$ it can be easily checked that 
\begin{equation}
\label{bound2fctG}
\int_0^{S(x)}|G_{\alpha}(t)|\, dt\leq \frac{3+2q+12 \alpha^2}{4(1-q)}\left(\frac{qx \sqrt{1-x^2}}{\sqrt{1-q x^2}}+S(x)\right)+\alpha(\alpha+1)\mathbf K(x,\sqrt{q})=g_{\alpha,q}(x),\quad x\in [0,1),
\end{equation} 
where ${\displaystyle \mathbf K(x,\sqrt{q})=\int_x^1 \frac{1}{\sqrt{(1-t^2)(1-q t^2)}}\, dt}.$ In particular, we have
\begin{equation}
\label{bound3fctG}
\int_0^{S(0)}|G_{\alpha}(t)|\, dt\leq \frac{3+2q+12 \alpha^2}{4(1-q)}\mathbf E(\sqrt{q})+
\alpha(\alpha+1)\mathbf K(\sqrt{q})= g_{\alpha,q}(0).
\end{equation}
Here, $\mathbf K(\cdot)$ and $\mathbf E(\cdot)$ are the Legendre Elliptic integrals of the first and the second kind,
given respectively by
$$\mathbf K(r)=\int_0^1 \frac{1}{\sqrt{(1-t^2)(1-r^2 t^2)}}\, dt,\qquad
\mathbf E(r)=\int_0^1 \sqrt{\frac{1-r^2 t^2}{1-t^2}}\, dt,\qquad r\in [0,1).$$
We have just proved the following Lemma.

\begin{lemma}
Under the above notations, consider two real numbers $ c>0,\, \alpha >-1$ and let $n$ be a positive 
integer so that ${\displaystyle q=\frac{c^2}{\chi_{n,\alpha}} <1.}$ If $V(\cdot)$ is  the function given  by
\eqref{Liouvile_transform2}, then it satisfies the differential equation 
\begin{equation}
\label{eq_diffV}
V"(s)+\left(\chi_{n,\alpha}+\frac{\frac{1}{4}-\alpha^2}{s^2}\right) V(s)= G_{\alpha}(s),\quad s\in (0, S(0)],
\end{equation} 
where $G_{\alpha}(\cdot)$ is given by \eqref{fctG} and satisfying  \eqref{bound2fctG} and \eqref{bound3fctG}.
\end{lemma}

Since $\ps$ has the same parity as $n,$ then it suffices to study the uniform approximation of the GPSWFs over the interval $[0,1].$ For this purpose, we use the following the weight and modulus  functions defined for any real $\alpha >-1$ as follows, see for example
[\cite{Olver}, p. 437],
\begin{equation}
\label{eight_fct}
 E_{\alpha}(x)=\left\{ \begin{array}{ll} (-Y_{\alpha}(x)/ J_{\alpha}(x))^{1/2} &\mbox{ if } 0<x\leq X_{\alpha}\\
 & \\ 1 &\mbox{ if } x\geq X_{\alpha},
 \end{array} \right.
\end{equation}
\begin{equation}
\label{modulus_fct}
M_{\alpha}(x)=\left\{ \begin{array}{ll} (2|Y_{\alpha}(x)| J_{\alpha}(x))^{1/2} &\mbox{ if } 0<x\leq X_{\alpha}\\
&\\
 (J^2_{\alpha}(x)+Y^2_{\alpha}(x))^{1/2} &\mbox{ if } x\geq X_{\alpha},  \end{array} \right.
\end{equation}
with $X_{\alpha}$ is the first zero of $J_{\alpha}(x)+Y_{\alpha}(x).$ The following proposition will be used 
for the error analysis study of the uniform approximation of the GPSWFs.

\begin{proposition} Under the above notation, for any real $\alpha \geq -\frac{1}{2},$ we have 
\begin{equation}
\label{bound_kernel}
\sup_{x >0} x M_{\alpha}^2(x) \leq m_{\alpha}= \left\{\begin{array}{cc}
2/\pi &\mbox{ if }  |\alpha |\leq \frac{1}{2} \\ & \\
\max\left( -2\alpha J_{\alpha}(\alpha) Y_{\alpha}(\alpha) +\frac{4\alpha}{\pi}; \alpha (J^2_{\alpha}(\alpha)+Y^2_{\alpha}(\alpha))\right), &\mbox{if } \alpha > \frac{1}{2}.\end{array} \right.
\end{equation}
\end{proposition}

\noindent
{\bf Proof:} We first note that from [\cite{Watson}, p.446], the function $ x \left(J^2_{\alpha}(x)+Y^2_{\alpha}(x)\right)$ is increasing if $|\alpha| \leq 1/2$ and it is decreasing if 
$\alpha > 1/2.$ Let   $|\alpha| \leq 1/2,$ since 
$$ 2 x |J_{\alpha}(x) Y_{\alpha}(x)| \leq x \left(J^2_{\alpha}(x)+Y^2_{\alpha}(x)\right),\quad x>0$$
and since $J_{\alpha}(x), Y_{\alpha}(x)$ have the asymptotic behaviours as $x\rightarrow +\infty,$ given by 
\begin{equation}
\label{asympt_infty}
J_{\alpha}(x) \sim \sqrt{\frac{2}{\pi x}}\cos \left(x-\alpha \frac{\pi}{2}-\frac{\pi}{4}\right),\quad 
Y_{\alpha}(x) \sim \sqrt{\frac{2}{\pi x}}\sin \left(x-\alpha \frac{\pi}{2}-\frac{\pi}{4}\right),
\end{equation}
then we have 
\begin{equation}\label{bound1_kernel}
\sup_{x >0} x M_{\alpha}^2(x) \leq \lim_{x\rightarrow +\infty} x \left(J^2_{\alpha}(x)+Y^2_{\alpha}(x)\right) =2/\pi,\quad |\alpha | \leq 1/2.
\end{equation}
Next, let  $\alpha > 1/2$ and let $j_{\alpha,1},  j'_{\alpha,1}$ and $y_{\alpha,1}$ denote the first zeros of 
$J_{\alpha}(x), J'_{\alpha}(x)$ and $Y_{\alpha}(x),$ respectively. It is known, see for example 
[\cite{Watson}, p. 487] that 
$$ \alpha < \sqrt{\alpha (\alpha+2)} < j'_{\alpha,1} < y_{\alpha,1}< j_{\alpha,1}.$$
Moreover, from the following asymptotic behaviours as $x\rightarrow 0^+,$ of the Bessel functions, given by
\begin{equation}
\label{asympt_zero}
J_{\alpha}(x) \sim \frac{1}{\Gamma(\alpha+1)} \left(\frac{x}{2}\right)^{\alpha},\quad 
Y_{\alpha}(x) = -\frac{1}{\pi}\Gamma(\alpha) \left(\frac{x}{2}\right)^{-\alpha},\quad x\rightarrow 0^+,
\end{equation}
one concludes that the function $ x J_{\alpha}(x) |Y_{\alpha}(x)|= - x J_{\alpha}(x) Y_{\alpha}(x) $ is positive and bounded over the interval $(0,\alpha].$  

Next, we check that for any $\alpha > 1/2,$ $ X_{\alpha} > \alpha ,$ where $X_{\alpha}$ is the first 
root  of $ J_{\alpha}(x)+J_{\alpha}(x)=0.$ To this end, we first  note that from 
[\cite{Andrews}, p.201],  the Wronskian of $J_{\alpha}, Y_{\alpha}$ is given by
\begin{equation}
\label{Wronskian1}
 W(J_{\alpha}, Y_{\alpha})(x)= J_{\alpha}(x) Y'_{\alpha}(x)-J'_{\alpha}(x)Y_{\alpha}(x)=\frac{2}{\pi x},\quad x>0.
\end{equation}
Consequently,  for any $\alpha > -1,$ we have
\begin{equation}
\label{deriv_YJ}
\frac{\partial }{\partial x}\left(\frac{-Y_{\alpha}(x)}{J_{\alpha}(x)}\right) = -\frac{2}{\pi x J^2_{\alpha}(x)} <0,\quad x >0.
\end{equation}
 Also, note that from [\cite{Olver}, p.438],  $X_{\alpha}$ has the following 
asymptotic formula, valid for large values of  the parameter $\alpha,$
$$ X_{\alpha} =\alpha + c (\alpha/2)^{1/3} +O(\alpha^{-1/3}),\quad c\approx 0.366.$$
Hence, there exists $\alpha_0 >0 ,$ so that $X_{\alpha} > \alpha,$ whenever $ \alpha \geq \alpha_0.$ Consequently, we have 
$$ - \frac{ Y_{\nu}(\nu)}{J_{\nu}(\nu)} \geq  - \frac{ Y_{\nu}(X_{\nu})}{J_{\nu}(X_{\nu})}=1,\quad\forall\, \nu \geq \alpha_0,$$
which means that ${\displaystyle \lim_{\nu\rightarrow +\infty}  - \frac{ Y_{\nu}(\nu)}{J_{\nu}(\nu)}\geq 1.}$ 
On the other hand, from [\cite{Watson}, p.515],  we have 
$$ \frac{\partial }{\partial \nu}\left(\frac{-Y_{\nu}(\nu)}{J_{\nu}(\nu)}\right) <0,\quad \nu >0.$$
Consequently, for any $\alpha \geq 1/2,$ we have 
$$ - \frac{ Y_{\alpha}(\alpha)}{J_{\alpha}(\alpha)} \geq \lim_{\nu\rightarrow +\infty}  - \frac{ Y_{\nu}(\nu)}{J_{\nu}(\nu)}\geq 1 = - \frac{ Y_{\alpha}(X_{\alpha})}{J_{\alpha}(X_{\alpha})},$$
which means that $ X_{\alpha}> \alpha,$ whenever  $\alpha \geq 1/2.$ Hence, for $0<x<\alpha,$  by integrating \eqref{deriv_YJ} over the interval $[x,\alpha]$ and using the fact that the function $x J^2_{\alpha}(x)$ is increasing,
one gets
\begin{eqnarray}
\label{bound1_JY}
-2 x J_{\alpha}(x) Y_{\alpha}(x)&=& 2 x J^2_{\alpha}(x) \frac{- Y_{\alpha}(\alpha)}{J_{\alpha}(\alpha)}+
\frac{4}{\pi} \int_x^{\alpha} \frac{x J^2_{\alpha}(x)}{t J^2_{\alpha}(t)}\, dt\nonumber\\
&\leq & 2\alpha J^2_{\alpha}(\alpha) \frac{- Y_{\alpha}(\alpha)}{J_{\alpha}(\alpha)}+ \frac{4\alpha}{\pi},\quad 0<x \leq \alpha.
\end{eqnarray} 
On the other hand, since $2 x |J_{\alpha}(x) Y_{\alpha}(x)\leq x \left(J^2_{\alpha}(x)+Y^2_{\alpha}(x)\right),$ since this later is decreasing for $\alpha \geq 1/2$ and since $ X_{\alpha}>\alpha,$ then we have 
\begin{equation}
\label{bound2_JY}
\max\left(\sup_{x\in [\alpha, X_{\alpha}]} -x J_{\alpha}(x) Y_{\alpha}(x); \sup_{x\geq X_{\alpha}}
x \left(J^2_{\alpha}(x)+Y^2_{\alpha}(x)\right)\right)\leq \alpha \left(J^2_{\alpha}(\alpha)+Y^2_{\alpha}(\alpha)\right).
\end{equation}
Finally, by combining \eqref{bound1_JY} and \eqref{bound2_JY}, one gets the desired bound \eqref{bound_kernel}.$\qquad \Box$

Once we have Lemma 1 and proposition 1, one can prove the following theorem that provides us with the uniform approximation over $[0,1]$ of the  $\psi^{(\alpha)}_{n,c}.$

\begin{theorem}\label{thm1}
 let  $ c>0, \alpha \geq -1/2 $ be two real numbers and let $n\in \mathbf N$ be such that 
 $q=c^2/\chi_{n,\alpha} <1$ and $(1-q)\sqrt{\chi_{n,\alpha}} \geq \pi (\frac{7}{4}+3\alpha^2) m_{\alpha}.$  Then under the previous notations, one can write
\begin{equation}\label{uniform0}
    \psi^{(\alpha)}_{n,c}(x)=  A_{\alpha}(q) \frac{(\chi_{n,\alpha})^{1/4}\sqrt{S(x)}J_{\alpha}(\sqrt{\chi_{n,\alpha}} S(x))}
{(1-x^2)^{1/4+\alpha/2}(1-q x^2)^{1/4}}+ \mathcal{E}_{n,\alpha}(x),\quad 0\leq x\leq 1,
\end{equation}
Here, $A_{\alpha}(q)$ is a normalization constant and 
\begin{equation}\label{bounds1}
|\mathcal{E}_{n,\alpha}(x)|\leq  \varepsilon_{n,\alpha} A_{\alpha}(q)
\frac{(1-x^2)^{1/4}}{(1-q x^2)^{3/4}}
\frac{\chi_{n,\alpha}^{1/4} \sqrt{S(x)}M_{\alpha}(\sqrt{\chi_{n,\alpha}}S(x))}{(1-x^2)^{\alpha/2}E_{\alpha}(\sqrt{\chi_{n,\alpha}}S(x))},
\end{equation}
where,
\begin{equation}\label{bounds2} 
\varepsilon_{n,\alpha} = \frac{1}{(1-q)\sqrt{\chi_{n,\alpha}}} \pi (e-1) (7/4 +3\alpha^2) m_{\alpha}.  \end{equation}
Here, $ m_{\alpha}$ is as given by  \eqref{bound_kernel}. 
\end{theorem}

\noindent 
{\bf Proof:} We first recall that for $x\in [0,1],$  ${\displaystyle \psi^{(\alpha)}_{n,c}(x)=\frac{V(S(x))}{(1-x^2)^{1/4+\alpha/2}(1-q x^2)^{1/4}},}$ where $V(\cdot)$ is a bounded 
solution on $[0,S(0)]$ of the differential equation \eqref{eq_diffV}. On the other hand, the general solution of 
this later is given by 
\begin{eqnarray}\label{expression_V}
V(s) &=& A_{\alpha}(q) V_1(\sqrt{\chi_{n,\alpha}} s)+ B_{\alpha}(q) V_2(\sqrt{\chi_{n,\alpha}} s)
+\int_0^s \frac{\sqrt{st\chi_{n,\alpha}}}{W(V_1(\sqrt{\chi_{n,\alpha}} \cdot),V_2(\sqrt{\chi_{n,\alpha}} \cdot))(t)}\cdot \nonumber\\
&&\qquad\qquad\Big( J_{\alpha}(\sqrt{\chi_{n,\alpha}} t) Y_{\alpha}(\sqrt{\chi_{n,\alpha}} s)-J_{\alpha}(\sqrt{\chi_{n,\alpha}} s) Y_{\alpha}(\sqrt{\chi_{n,\alpha}} t)\Big) G_{\alpha}(t) V(t)\, dt.
\end{eqnarray}
Here, $V_1(t) = \sqrt{t} J_{\alpha}(t),\, V_2(t)= \sqrt{t} Y_{\alpha}(t).$ For the homogeneous solutions 
$V_1(\sqrt{\chi_{n,\alpha}} s), V_2(\sqrt{\chi_{n,\alpha}} s)$ of \eqref{eq_diffV}, one may refer to [\cite{Andrews}, p. 201]. Note that from \eqref{equS} and the asymptotic behaviours of $J_{\alpha}, Y_{\alpha},$ given by \eqref{asympt_zero}, one concludes that the function ${\displaystyle \frac{V_1(S(x))}{(1-x^2)^{1/4+\alpha/2} (1-q x^2)^{1/4}}}$ is bounded at $x=+1,$  or $s=S(1)=0,$ which is not the case for the function ${\displaystyle \frac{V_2(S(x))}{(1-x^2)^{1/4+\alpha/2} (1-q x^2)^{1/4}}}.$ Hence, in \eqref{expression_V},  we have $B_{\alpha}(q)=0.$
Moreover, from the expression of the Wronskian given by \eqref{Wronskian1}, one can easily check that
\begin{eqnarray}\label{expression2_V}
V(S(x))&=&  A_{\alpha}(q) {(\chi_{n,\alpha})^{1/4}\sqrt{S(x)}J_{\alpha}(\sqrt{\chi_{n,\alpha}} S(x))}
+\frac{\pi \sqrt{S(x)} }{2 \sqrt{\chi_{n,\alpha}}}\nonumber \\ 
&&\qquad\cdot \int_0^{S(x)} \sqrt{t}\Big( J_{\alpha}(\sqrt{\chi_{n,\alpha}} t) Y_{\alpha}(\sqrt{\chi_{n,\alpha}} S(x))-J_{\alpha}(\sqrt{\chi_{n,\alpha}} S(x)) Y_{\alpha}(\sqrt{\chi_{n,\alpha}} t)\Big) G_{\alpha}(t) V(t)\, dt\nonumber\\
&=& A_{\alpha}(q) {(\chi_{n,\alpha})^{1/4}\sqrt{S(x)}J_{\alpha}(\sqrt{\chi_{n,\alpha}} S(x))}+ R_{n,\alpha}(x).
\end{eqnarray} 
On the other hand, from [\cite{Olver}, p.450], a bound of  the reminder term $R_{n,\alpha}(x)$ is given as follows,
\begin{equation}
\label{Olver_bound}
|R_{n,\alpha}(x)| \leq A_{\alpha}(q)\chi_{n,\alpha}^{1/4}\sqrt{S(x)}\frac{M_{\alpha}(\sqrt{\chi_{n,\alpha}}S(x))}{E_{\alpha}(\sqrt{\chi_{n,\alpha}}S(x))}
\left(e^{\gamma_n(x)} -1\right),
\end{equation}
where ${\displaystyle \gamma_n(x)= \frac{\pi}{2 \sqrt{\chi_{n,\alpha}}} \int_0^{S(x)} t \sqrt{\chi_{n,\alpha}} M^2_{\alpha}(\sqrt{\chi_{n,\alpha}} t)
|G_{\alpha}(t)|\, dt.}$ Moreover, from \eqref{bound2fctG} and\eqref{bound_kernel}, we have 
\begin{equation}
\label{bound1gamma}
\gamma_n(x) \leq  \frac{\pi}{2 \sqrt{\chi_{n,\alpha}}} m_{\alpha} g_{\alpha,q}(x),\quad x\in [0,1].
\end{equation}
Also, since $${\displaystyle \mathbf K(x,\sqrt{q})\leq \frac{1}{\sqrt{1-q x^2}} \int_x^1\frac{dt}{\sqrt{1-t^2}}\leq
\frac{2}{1-q}\frac{\sqrt{1-x^2}}{\sqrt{1-q x^2}}}$$ and since from \cite{Bonami-Karoui2}, we  have 
${\displaystyle \frac{qx \sqrt{1-x^2}}{\sqrt{1-q x^2}}+S(x)\leq 2 \frac{\sqrt{1-x^2}}{\sqrt{1-q x^2}},}$ one gets 
$$|g_{\alpha,q}(x)|\leq \frac{\sqrt{1-x^2}}{\sqrt{1-q x^2}}\left(\frac{3+12\alpha^2+2q}{2(1-q)}+\frac{4(1-q)}{2(1-q)}\right).$$
Consequently, we have
\begin{equation}\label{boundg}
\frac{|g_{\alpha,q}(x)|}{(1-x^2)^{1/4}(1-q x^2)^{1/4}}\leq \frac{1}{1-q}\frac{(1-x^2)^{1/4}}{(1-qx^2)^{1/4}}(7/2 + 6\alpha^2).
\end{equation}
Moreover, since $(1-q)\sqrt{\chi_{n,\alpha}} \geq \pi (\frac{7}{4}+3\alpha^2) m_{\alpha},$
then for $x\in [0,1],$ we have $\gamma_n(x)\leq \gamma_n(0)\leq 1.$ Hence, we have 
\begin{equation}\label{bound2gamma}
\left(e^{\gamma_n(x)} -1\right)\leq \gamma_n(x) \frac{e^{\gamma_n(0)}-1}{\gamma_n(0)}\leq \gamma_n(x) (e-1),
\end{equation}
Finally, since ${\displaystyle \psi_{n,c}^{(\alpha)}(x)=(1-x^2)^{-\alpha/2-/4}(1-q x^2)^{-1/4} V(S(x)),}$ then by using 
\eqref{expression2_V}, \eqref{bound1gamma}, \eqref{boundg} and \eqref{bound2gamma}, one gets the desired 
result \eqref{uniform0}.$\qquad \Box$

\begin{remark}
Since $${\displaystyle \frac{M_{\alpha}(s)}{E_{\alpha}(s)}=\left\{\begin{array}{ll} \sqrt{2} J_{\alpha}(s) &\mbox{ if } 0<s\leq X_{\alpha}\\
&\\
\left(J^2_{\alpha}(s)+Y^2_{\alpha}(s)\right)^{1/2} &\mbox{ if }  s \geq X_{\alpha},\end{array}\right.}$$ then from the 
asymptotic behaviour of $J_{\alpha}(x),$ given by \eqref{asympt_zero}, one concludes that the quantity in the reminder term \eqref{bounds1},  given by 
${\displaystyle \frac{\chi_{n,\alpha}^{1/4} \sqrt{S(x)}M_{\alpha}(\sqrt{\chi_{n,\alpha}}S(x))}{(1-x^2)^{\alpha/2}E_{\alpha}(\sqrt{\chi_{n,\alpha}}S(x))}}$ 
is bounded on $[0,1).$ 
\end{remark}

Note that since $\varepsilon_n(1)=0,$ and since from \cite{Bonami-Karoui2}, we have 
${\lim_{x\rightarrow 1} S(x)/\sqrt{(1-x^2)(1-q x^2)} =1},$ then by using the asymptotic behaviour of 
$J_{\alpha}(\sqrt{\chi_{n,\alpha}} S(x))$ as $x\rightarrow 1^{-},$ which is obtained from 
\eqref{asympt_zero}, one concludes that the normalisation constant $A_{\alpha}(q)$ of Theorem 1, is given by
\begin{equation}
\label{constantA}
A_{\alpha}(q)= \frac{2^{\alpha} \Gamma(1+\alpha)}{(1-q)^{\alpha/2} \chi_{n,\alpha}^{1/4+\alpha/2}} 
\psi_{n,c}^{\alpha}(1).
\end{equation}

Next, we give an accurate explicit approximation of the normalisation constant $A_{\alpha}(q),$ so that 
the  $\psi_{n,c}^{\alpha}$ are normalized by the requirement that 

\begin{equation}
\label{normalpsi}
\| \psi^{(\alpha)}_{n,c}\|^2_{L^2([-1,1],\omega_{\alpha})}=\int_{-1}^1  \left(\psi^{(\alpha)}_{n,c}(x)\right)^2 \, \omega_{\alpha}(x)\, dx =1,\qquad \omega_{\alpha}(x)=(1-x^2)^{\alpha}.
\end{equation}
To this end, we first define the following two constants depending on $\alpha,$
\begin{equation}
\label{constants}
\mu_{\alpha}=|\alpha^2-\frac{1}{4}|, \quad       c_{\alpha} =\left\{\begin{array}{ll}
\sqrt{2/\pi} &\mbox{ if } |\alpha|\leq 1/2\\
0.675\sqrt{\alpha^{1/3}+\frac{1.9}{\alpha^{1/3}}+
\frac{1.1}{\alpha}}&\mbox{ if } \alpha >1/2.\end{array}\right.
\end{equation}
The following lemma is essential in the estimate of the normalization constant $A_{\alpha}(q).$

\begin{lemma}
Let $ \alpha \geq  -1/2, $ then  for any $x>0,$ we have 
\begin{equation}
\label{Primitive1}
\int_{0}^{x} t J^2_{\alpha}(t)\, dt = \frac{x^2}{2}\left[ J^2_{\alpha}(x)+J^2_{\alpha+1}(x) -\frac{2\alpha}{x}J_{\alpha}(x)J_{\alpha+1}(x) \right] = \frac{x}{\pi}+\eta_{\alpha}(x),
\end{equation}
where 
\begin{equation}
\label{boundeta}
\sup_{x\geq 0} |\eta_{\alpha}(x)| \leq M_{\alpha}=\max\left(\frac{1}{\pi}, c_{\alpha}^2-\frac{1}{\pi}, \kappa_{\alpha}\right)
\end{equation}
with $${\kappa_{\alpha}= \frac{4}{5}\sqrt{\frac{2}{\pi}} (\mu_{\alpha}+\mu_{\alpha+1})+\frac{8}{25}(\mu^2_{\alpha}+\mu^2_{\alpha+1})+|\alpha| c_{\alpha} c_{\alpha+1}}.$$
\end{lemma}

\noindent
{\bf Proof:} The  first equality in \eqref{Primitive1} is a consequence of the following 
identity, see [\cite{NIST}, p.241]
$$\int_0^x t J^2_{\alpha}(t)\, dt = \frac{x^2}{2}\left[ J^2_{\alpha}(x)-J_{\alpha-1}(x)J_{\alpha+1}(x)\right],\quad \alpha > -1/2,$$ combined with the well known identity 
$$J_{\alpha-1}(x)=\frac{2\alpha}{x}J_{\alpha}(x)-J_{\alpha+1}(x).$$ Moreover, it has been shown in \cite{Krasikov}, that for $\alpha \geq -1/2,$ we have
\begin{equation}
\label{bound0J}
\sup_{x\geq 0} x^{3/2}\left|J_{\alpha}(x)-\sqrt{\frac{2}{\pi x}} \Big[ \cos\left(x-(\alpha+1/2)\frac{\pi}{2}\right)\right|\leq 
\frac{4}{5} \mu_{\alpha}.
\end{equation}
Hence, by using the previous inequality, one gets 
\begin{equation}
\label{bound1J}
\left| J_{\nu}^2(x)-\frac{2}{\pi x} \cos\big(x-(\nu+\frac{1}{2})\frac{\pi}{2}\big)\right|\leq \frac{4}{5 x^{3/2}}\mu_{\alpha}
\left( |J_{\nu}(x)|+\sqrt{\frac{2}{\pi x}}\right),\quad x\geq 1,\quad \nu = \alpha,\, \alpha+1.
\end{equation}
Moreover, from \eqref{bound0J}, one gets 
\begin{equation}
\label{bound2J}
|J_{\nu}(x)| \leq \frac{4}{5} \mu_{\nu}+\sqrt{\frac{2}{\pi x}},\quad x\geq 1.
\end{equation}
By using the previous two inequalities, one obtains
$$\left| J^2_{\alpha}(x)+J^2_{\alpha+1}(x)-\frac{2}{\pi x}\right|\leq \frac{8}{5} \sqrt{\frac{2}{\pi}} \frac{1}{x^2}(\mu_{\alpha}+\mu_{\alpha+1})+\frac{16}{25 x^3}(\mu^2_{\alpha}+\mu^2_{\alpha+1}),\quad x\geq 1.$$
Hence, we have
\begin{eqnarray}
\label{Ineq2.1}
\lefteqn{\frac{x^2}{2}\left| J^2_{\alpha}(x)+J^2_{\alpha+1}(x)-\frac{2\alpha}{x} J_{\alpha}(x)J_{\alpha+1}(x)-\frac{2}{\pi x}\right|}\nonumber\\
&&\qquad\qquad \leq \frac{4}{5} \sqrt{\frac{2}{\pi}} \frac{1}{x^2}(\mu_{\alpha}+\mu_{\alpha+1})+\frac{8}{25 x}(\mu^2_{\alpha}+\mu^2_{\alpha+1})+|\alpha| x|J_{\alpha}(x)J_{\alpha+1}(x)| ,\quad x\geq 1.
\end{eqnarray}
On the other hand, it has been shown in \cite{Olenko}, that 
\begin{equation}
\label{boundJ}
\sup_{x\geq 0} \sqrt{x} |J_{\alpha}(x)| \leq c_{\alpha},
\end{equation}
where an upper bound of $c_{\alpha}$ is given by \eqref{constants}.
Finally, by combining the previous two inequalities, one gets a bound for $|\eta_{\alpha}(x)|$ for $x\geq 1.$ To get 
a bound 
$\eta_{\alpha}(x)$ over the interval $[0,1],$ it suffices to note that from \eqref{Primitive1}, we have
$$\sup_{x\in [0,1]} |\eta'_{\alpha}(x)|=\sup_{x\in [0,1]} \left| x J^2_{\alpha}(x) -\frac{1}{\pi}\right|\leq\max\left(\frac{1}{\pi}, c^2_{\alpha}-\frac{1}{\pi}\right).$$
Since $\eta_{\alpha}(0)=0,$ then the previous bound is also valid for 
${\displaystyle \sup_{x\in [0,1]} |\eta_{\alpha}(x)| },$ that is 
\begin{equation}
\label{bound3J}
|\eta_{\alpha}(x)|\leq \max\left(\frac{1}{\pi}, c^2_{\alpha}-\frac{1}{\pi}\right),\quad 0\leq x\leq 1.
\end{equation}
Finally, to conclude for the proof of the lemma, it suffices to combine 
\eqref{Ineq2.1}, \eqref{boundJ} and \eqref{bound3J}.$\quad\Box$

The following lemma provides us with an explicit estimate of the weighted $L^2([0,1],\omega_{\alpha})-$norm of 
$\widetilde \psi^{(\alpha)}_{n,c},$ the uniform approximation of the GPSWFs, given in Theorem 1, by
\begin{equation}
\label{approxpsi}
\widetilde \psi^{(\alpha)}_{n,c}(x)= A_{\alpha}(q) \frac{(\chi_{n,\alpha})^{1/4}\sqrt{S(x)}J_{\alpha}(\sqrt{\chi_{n,\alpha}} S(x))}
{(1-x^2)^{1/4+\alpha/2}(1-q x^2)^{1/4}},\quad x\in [0,1).
\end{equation}

\begin{lemma}
Under the previous notations, let $c>0,$ $\alpha \geq -1/2$ be two real numbers. Then, for any  $n\in \mathbf N$
with $q=c^2/\chi_{n,\alpha}<1,$ we have 
\begin{equation}
\label{norm2tpsi}
\left| \|\widetilde \psi^{(\alpha)}_{n,c}\|^2_{L^2([0,1],\omega_{\alpha})}- A^2_{\alpha}(q)
\frac{\mathbf K (\sqrt{q})}{\pi}\right|\leq A^2_{\alpha}(q)\frac{M_{\alpha}}{(1-q) \sqrt{\chi_{n,\alpha}}},
\end{equation} 
where $M_{\alpha}$ is given by \eqref{boundeta}. 
\end{lemma}

\noindent
{\bf Proof: } We first write $\|\widetilde \psi^{(\alpha)}_{n,c}\|^2_{L^2([0,1],\omega_{\alpha})}$ as follows

\begin{eqnarray}\label{Eq2.1}
\int_0^1 \left(\widetilde \psi^{(\alpha)}_{n,c}\right)^2(x) \omega_{\alpha}(x)\, dx &=& A^2_{\alpha}(q)
\int_0^1 \sqrt{\chi_{n,\alpha}} S(x) J^2_{\alpha}(\sqrt{\chi_{n,\alpha}} S(x)) S'(x) \frac{1}{1-q x^2}\, dx\nonumber\\
&=& A^2_{\alpha}(q) \int_0^1 F'_n(x) \frac{1}{1-q x^2} \, dx,
\end{eqnarray}
with 
$$ F_n(x)= -\int_x^1 \sqrt{\chi_{n,\alpha} } S(t) J^2_{\alpha}(\sqrt{\chi_{n,\alpha}} S(t)) S'(t)\, dt.$$
Since $S(1)=0,$ then a change of variable and Lemma 2, give us
$$ F_n(x)= \frac{1}{\sqrt{\chi_{n,\alpha}}}\int_0^{\sqrt{\chi_{n,\alpha}}S(x)} u J^2_{\alpha}(u)\, du = \frac{S(x)}{\pi}+
\frac{\eta_{\alpha}(\sqrt{\chi_{n,\alpha}} S(x))}{\sqrt{\chi_{n,\alpha}}}.$$
Since $F_n(1)=0,$ then by using the previous equality and  integrations by parts applied to the integral in \eqref{Eq2.1}, one gets 
\begin{eqnarray}
\int_0^1 F'_n(x) \frac{1}{1-q x^2} \, dx&=& \frac{1}{\pi} \int_0^1 S(x) \frac{2q x}{(1-q x^2)^2}\, dx +\frac{1}{\sqrt{\chi_{n,\alpha}}} \int_0^1 \eta_{\alpha}(\sqrt{\chi_{n,\alpha}} S(x))\frac{2q x}{(1-q x^2)^2}\, dx\nonumber\\
&=& -F_n(0)+\frac{1}{\pi} \int_0^1 \sqrt{\frac{1-q x^2}{1-x^2}} \frac{1}{1-q x^2}\, dx + \frac{1}{\sqrt{\chi_{n,\alpha}}} I_{\alpha,q}\nonumber\\
&=& \frac{\mathbf K(\sqrt{q})}{\pi} + \frac{1}{\sqrt{\chi_{n,\alpha}}} I_{\alpha,q},
\end{eqnarray}
where ${\displaystyle I_{\alpha,q}=-\frac{\eta_{\alpha}(\sqrt{\chi_{n,\alpha}} S(0))}{\sqrt{\chi_{n,\alpha}}}+\frac{1}{\sqrt{\chi_{n,\alpha}}}\int_0^1 \eta_{\alpha}(\sqrt{\chi_{n,\alpha}} S(x))\frac{2q x}{(1-q x^2)^2}\, dx.}$
Moreover, from  \eqref{boundeta}, one gets ${\displaystyle |I_{\alpha,q}|\leq M_{\alpha} \frac{1}{1-q}.}$
This concludes the proof of the lemma.$\quad \Box$

The following proposition provides us with an estimate of the normalisation constant $A_{\alpha}(q)$ under the condition
\eqref{normalpsi}.

\begin{proposition} Let $c>0$ and $\alpha \geq -1/2,$  then there exists $N_{\alpha}\in  \mathbb N,$
such that for any $n\geq N_{\alpha},$ we have $q=\frac{c^2}{\chi_{n,\alpha}}\leq q_0<1,$ $(1-q)\sqrt{\chi_{n,\alpha}} \geq \pi (\frac{7}{4}+3\alpha^2) m_{\alpha}.$ Moreover, there exists a constant 
$C_{\alpha}$ depending only on $\alpha$ and   such that  for $n\geq N_{\alpha},$ we have
\begin{equation}\label{boundsA}
\sqrt{\frac{\pi}{2\mathbf K(\sqrt{q})}}\frac{1}{1+\epsilon_{n,\alpha}C_{\alpha}}\leq A_{\alpha}(q) \leq \sqrt{\frac{\pi}{2\mathbf K(\sqrt{q})}}\frac{1}{1-\epsilon_{n,\alpha} C_{\alpha}}, 
\end{equation}
where $\epsilon_{n,\alpha}$ is given by \eqref{bounds2}.
\end{proposition}

\noindent
{\bf Proof:} From \eqref{bounds1} and the expression of ${\displaystyle \frac{\sqrt{t} M_{\alpha}(t)}{E_{\alpha}(t)},}$
one can easily check that 
\begin{eqnarray*}
\|\mathcal E_{n,\alpha}\|_{L^2([0,1],\omega_{\alpha})}&\leq & \epsilon_{n,\alpha} \max(\sqrt{2} c_{\alpha},\sqrt{m_{\alpha}}) A_{\alpha}(q) \int_0^1 \frac{(1-x^2)^{1/4}}{(1-q x^2)^{3/4}}\, dx\\
&\leq & \epsilon_{n,\alpha} \max(\sqrt{2} c_{\alpha},\sqrt{m_{\alpha}}) A_{\alpha}(q) \frac{\pi}{2}.
\end{eqnarray*}
Since, $\|\ps\|_{L^2([0,1],\omega_{\alpha})}=\frac{1}{\sqrt{2}},$ then the previous inequality implies
\begin{equation}\label{normal2}
\left| \|\widetilde \psi^{(\alpha)}_{n,c}\|_{L^2([0,1],\omega_{\alpha})}- \frac{1}{\sqrt{2}}\right|\leq
\epsilon_{n,\alpha} C_{\alpha} A_{\alpha}(q).
\end{equation}
Moreover, since $\mathbf K(\sqrt{q}) \geq \frac{\pi}{2}$ and since $A_{\alpha}(q) >0,$ then $ \left|\|\widetilde \psi^{(\alpha)}_{n,c}\|_{L^2([0,1],\omega_{\alpha})}+ A_{\alpha}(q)
\frac{\mathbf K (\sqrt{q})}{\pi}\right|\geq \frac{A_{\alpha}(q)}{\sqrt{2}}.$ By combining this last inequality 
with \eqref{norm2tpsi}, one gets
\begin{equation}\label{normal3}
\left| A_{\alpha}(q)\sqrt{\frac{\mathbf K(\sqrt{q})}{\pi}}- \frac{1}{\sqrt{2}}\right|\leq
\epsilon_{n,\alpha} C_{\alpha} A_{\alpha}(q).
\end{equation}
Finally, since $\sqrt{\frac{\pi}{\mathbf K(\sqrt{q})}}\leq \sqrt{2},$ then the previous inequality gives us
the desired inequalities \eqref{boundsA}.$\qquad \Box$

\begin{remark}
By combining the results of Theorem 1 and Proposition 2, one obtains the following uniform asymptotic approximation
of $\ps$ in terms of the Bessel function $J_{\alpha}(\cdot),$
$$\ps = \sqrt{\frac{\pi}{2\mathbf K(\sqrt{q})}} \frac{(\chi_{n,\alpha})^{1/4}\sqrt{S(x)}J_{\alpha}(\sqrt{\chi_{n,\alpha}} S(x))}
{(1-x^2)^{1/4+\alpha/2}(1-q x^2)^{1/4}}+ \mathbf E_{n,\alpha}(x), \quad 0\leq x\leq 1,$$
where
$$\mathbf E_{n,\alpha}(x) = \left(A_{\alpha}(q)-\sqrt{\frac{\pi}{2\mathbf K(\sqrt{q})}}\right)\frac{(\chi_{n,\alpha})^{1/4}\sqrt{S(x)}J_{\alpha}(\sqrt{\chi_{n,\alpha}} S(x))}
{(1-x^2)^{1/4+\alpha/2}(1-q x^2)^{1/4}}+ \mathcal{E}_{n,\alpha}(x).$$
Here, $A_{\alpha}(q)$ and $\mathcal{E}_{n,\alpha}(x)$ are as given by Theorem 1.
\end{remark}

\section{Uniform approximation of the eigenfunctions in terms of Jacobi polynomials.} 

In this paragraph, we show that for a given real number $c>0$ and any real $0< \alpha < 3/2,$ the GPSWFs 
$\psi^{(\alpha)}_{n,c}$ are uniformly  approximated by the normalized Jacobi polynomial $\widetilde P_n^{(\alpha,\alpha)}.$ For this purpose, we first need the following mathematical preliminaries on Jacobi polynomials and Jacobi functions of the second kind.

 \subsection{Preliminaries on Jacobi polynomials and  Jacobi functions of the second kind.}
 
We recall that for two real numbers $\alpha, \beta >-1,$  the Jacobi polynomials $\J_k$ are given by the recurrence 
formula
$$
\J_{k+1}(x)= (A_k x + B_k)\J_{k}(x) -C_k \J_{k-1}(x),\quad x\in [-1,1],
$$
where  $\J_0(x)=1,\quad \J_1(x)=\frac{1}{2}(\alpha+\beta+2)x +\frac{1}{2}(\alpha+\beta)$ and where 
$
A_k =\frac{(2k+\alpha+\beta+1)(2k+\alpha+\beta+2)}{2 (k+1)(k+\alpha+\beta+1)},$ $B_k=\frac{(\alpha^2-\beta^2)(2k+\alpha+\beta+1)}{2(k+1)(k+\alpha+\beta+1)(2k+\alpha+\beta)},\,\, 
C_k=\frac{(k+\alpha)(k+\beta)(2k+\alpha+\beta+2)}{(k+1)(k+\alpha+\beta+1)(2k+\alpha+\beta)}.$
The normalized Jacobi polynomial of degree $k,$ denoted by $\wJ_k$ and satisfying the condition   
${\displaystyle \int_{-1}^1 (\wJ_k(y))^2 (1-y)^{\alpha}(1+y)^{\beta}\, dy =1}$ are given by 
\begin{equation}\label{JacobiP}
\wJ_{k}(x)= \frac{1}{\sqrt{h_k^{(\alpha,\beta)}}}\J_k(x),\quad h_k^{(\alpha,\beta)}=\frac{2^{\alpha+\beta+1}\Gamma(k+\alpha+1)\Gamma(k+\beta+1)}{k!(2k+\alpha+\beta+1)\Gamma(k+\alpha+\beta+1)}.
\end{equation}
Here, $\Gamma(\cdot)$ is the gamma function that satisfies  
the following useful inequalities, see \cite{Batir} 
\begin{equation}\label{Ineq2}
\sqrt{2e} \left(\frac{x+1/2}{e}\right)^{x+1/2}\leq \Gamma(x+1)\leq \sqrt{2\pi} \left(\frac{x+1/2}{e}\right)^{x+1/2},\quad x>0. 
\end{equation}
Note that $P_n^{(\alpha,\beta)}$ is the bounded solution of the following second order differential equation,
$$(1-x^2) y"(x)+(\beta-\alpha-(\alpha+\beta+2)x)y'(x)+n(n+\alpha+\beta+1)y(x)=0,\quad x\in (-1,1).$$
A second linearly independent solution of the previous differential equation is given by the Jacobi function of the second kind, denoted by $Q_n^{(\alpha,\beta)}$ and defined by
$$Q_n^{(\alpha,\beta)}(x)=Q_0^{(\alpha,\beta)}(x)P_n^{(\alpha,\beta)}(x)-\frac{W_{n-1}^{(\alpha,\beta)}(x)}{(1-x)^{\alpha}(1+x)^{\beta}},$$
where
$$Q_0^{(\alpha,\beta)}(x)= \int_0^x \frac{(1+\alpha+\beta) dt }{(1-t)^{1+\alpha}(1+t)^{1+\beta}}+\Lambda_{\alpha \beta},$$ $$ W_{n-1}^{(\alpha,\beta)}(x)=\frac{\Gamma(\alpha+\beta+2)}{2^{\alpha+\beta+1}\Gamma(\alpha+1)\Gamma(\beta+1)}\int_{-1}^1
\frac{P_n^{(\alpha,\beta)}(x)-P_n^{(\alpha,\beta)}(s)}{x-s}(1-s)^{\alpha}(1+s)^{\beta}\, ds.  $$
Here, $\Lambda_{\alpha \beta}$ is a constant depending on $\alpha,\,\beta$ and $W_{n-1}^{(\alpha,\beta)}$ is the first associated polynomial. For more details, see \cite{Grosjean}. 
It is interesting to note that  the Jacobi polynomial and the Jacobi function of the second kind $Q_n^{(\alpha,\beta)}$
satisfy the following local estimates, see \cite{Nevai}
\begin{eqnarray}\label{estimatesJacobi}
|Q_n^{(\alpha,\beta)}(x)| &\leq& \frac{C}{(1-x)^{\alpha}(1+x)^{\beta}} \left(\sqrt{1-x}+n^{-1}\right)^{\alpha-1/2} \left(\sqrt{1+x}+n^{-1}\right)^{\beta-1/2},\quad -1<x<1,\nonumber\\
|P_n^{(\alpha,\beta)}(x)| &\leq& C \left(\sqrt{1-x}+n^{-1}\right)^{-\alpha-1/2} \left(\sqrt{1+x}+n^{-1}\right)^{-\beta-1/2},\quad -1< x< 1,
\end{eqnarray}
where $C$ is a fixed constant, not depending on the parameters $n,\,\alpha,\, \beta.$
 Note that if $x_{n,k}$ are the $n$ zeros of $P_n^{(\alpha,\beta)},$ arranged in the decreasing order $-1 < x_{n,n}<\cdots < x_{n,1}< 1,$ then it has been shown in \cite{Area} that if $\alpha, \beta > -1/2,$ then $Q_n^{(\alpha,\beta)}$ has $n+1$ zeros in $(-1,1),$ denoted
by $t_{n,k},$ arranged in the decreasing order and  satisfying the following interlacing property
\begin{equation}\label{interlacing}
x_{n,k+1} < t_{n,k} < x_{n,k},\quad k=1,\ldots,n-1,\qquad t_{n,0}\in (x_{n,1}, 1).
\end{equation}
Also, from [\cite{Szego}, p. 192], an asymptotic formula for the zeros $x_{n,k}$ is given by
\begin{equation}\label{asymp_xk}
x_{n,k}=\cos \theta_{n,k},\quad \lim_{n\rightarrow +\infty} n \theta_{n,k}=j_{k,\alpha},
\end{equation}
where $j_{k,\alpha}$ is the $k-$th positive zero of the Bessel function $J_{\alpha}(\cdot).$ 
Moreover, from [\cite{Watson}, p. 506], for fixed $\alpha >-1$ and for large enough integer $k,$ we have 
the following asymptotic approximation of  $j_{k,\alpha},$
\begin{equation}\label{asymp_jk}
j_{k,\alpha}= k\pi +\frac{\pi}{2}\left(\alpha-\frac{1}{2}\right)-\frac{4\alpha^2-1}{8(k\pi+\frac{\pi}{2}(\alpha-\frac{1}{2}))}+ O(k^{-3}).
\end{equation}

\subsection{Uniform approximation in terms of Jacobi polynomials}
In the sequel, we let $C_{\alpha}(q_0)$ denote a generic constant depending on 
$\alpha$ and $0<q_0<1.$  The following theorem provides us with the approximation of $\ps$ by $\pn.$

\begin{theorem}
Let $c>0$ and $0<\alpha < 3/2,$  then there exists $N_{\alpha}\in  \mathbb N,$
such that for any $n\geq N_{\alpha},$ we have $q=\frac{c^2}{\chi_{n,\alpha}}\leq q_0<1,$ $(1-q)\sqrt{\chi_{n,\alpha}} \geq \pi (\frac{7}{4}+3\alpha^2) m_{\alpha}.$ Moreover, there exists a constant 
$C_{\alpha}(q_0)$ depending only on $\alpha$ and $q_0$ and  such that  for $n\geq N_{\alpha},$ we have
\begin{equation}\label{approx1P}
 \left|\ps(x)-A_n \pn (x)\right|\leq C_{\alpha}(q_0) \frac{  c^2}{n+2\alpha+1},\quad \forall\,\, x\in [-1,1],
\end{equation}
where  $A_n$ is the normalization constant, satisfying
\begin{equation}\label{approx2P}
|1- A_n | \leq C_{\alpha}(q_0) \frac{c^2}{2n+2\alpha+1}.
\end{equation}
\end{theorem}

\noindent
{\bf Proof:}  We will only prove the previous approximation result on $[0,1],$ since the same proof is used
on the interval $[-1,0].$ We  first rewrite the differential equation
 governing $\ps$ as follows
\begin{equation}\label{Eq1}
(1-x^2) \psi"(x) - 2 x \psi'(x) +\chi_n(0) \psi(x) = \big(\chi_{n,\alpha}(0)-\chi_{n,\alpha}(c) +c^2 x^2\big) \psi(x),\quad x\in [0,1],
\end{equation}
where $\chi_{n,\alpha}(0)=n(n+2\alpha+1).$ Note that the homogeneous equation associated with the previous differential equation has $\pn$ and $\qn$ as the two linearly
independent solutions. 
By the method of variation of constants, the bounded solution $\ps$ of the previous equation is written as
\begin{equation}\label{Eq2}
\ps(x)= A_n \pn+ B_n \qn +\int_x^1 \frac{k_{n}(x,y) G(y) \ps(y)}{W(\pn,\qn)(y)}\, dy=
A \pn(x)+ B \qn(x) + R_{n,\alpha}(x). 
\end{equation}
where $A_n, B_n$ are constants and 
\begin{equation}\label{Eq3}
k_{n}(x,y)= \pn(x)\qn(y)- \pn(y)\qn(x).
\end{equation}
Also, since $
G(y)= \chi_{n,\alpha}(0)-\chi_{n,\alpha}(c) + c^2 y^2$ and since 
$-c^2 y^2\leq \chi_{n,\alpha}(0)-\chi_{n,\alpha}(c)\leq 0,$ then we have 
\begin{equation}\label{Eq4}
|G(y)|\leq c^2 y^2,\quad y\in [0,1].
\end{equation}
Also, from [\cite{Erdelyi}, p.171] and taking into account the normalization constant 
$h_n^{(\alpha,\alpha)},$ given by \eqref{JacobiP},  as well as the bounds of gamma function, given by 
\eqref{Ineq2},   one gets the following estimate of the Wronskian $W(\pn,\qn)(y)$ 
\begin{eqnarray}\label{Wronskian}
|W(\pn,\qn)(y)|&=& \frac{2^{2\alpha}}{h_n^{(\alpha,\alpha)}} \frac{\Gamma^2(n+\alpha+1)}{\Gamma(n+1)\Gamma(n+2\alpha+1)}\frac{1}{(1-y^2)^{1+\alpha}}\nonumber\\
&\geq & C_{\alpha} \frac{2n+2\alpha+1}{(1-y^2)^{1+\alpha}}.
\end{eqnarray}

  Next, we prove that the kernel $K_n(x,y)= (1-y^2)^{1+\alpha} k_n(x,y)$
   is bounded on the set $\{x,y \in [0,1];\,\, y\geq x\}.$ For this purpose, we first note that from
the interlacing property of zeros of $\pn$ and $\qn,$ given by \eqref{interlacing}, as well as from
the asymptotic zeros locations of Jacobi polynomials, given in the previous paragraph,   one concludes that there exists a constant 
$\gamma>0$ and a positive integer $N_{\alpha}\in \mathbb N$ such that 
$$ 1- \frac{\gamma}{n} \leq t_{n,0} < 1,\quad \forall n\geq N_{\alpha}.$$
Recall that $t_{n,0}$ is the largest zero of $\qn$ in $(0,1).$ On the other hand, from [\cite{Szego}, p.67], the function ${\displaystyle u_{\alpha}(x)= (1-x^2)^{(1+\alpha)/2} \qn(x)}$ is a solution of the following differential equation
$$u_{\alpha}" (x) + g_{n,\alpha}(x) u_{\alpha}(x)=0,\quad g_{n,\alpha}(x)=\frac{1}{4}\left( \frac{(1-\alpha^2)(1+x^2)}{(1-x^2)^2}+\frac{4n(n+2\alpha+1)+2(1+\alpha)^2}{1-x^2}\right).$$
Since
$$ g'_{n,\alpha}(x)=\frac{8 x}{(1-x^2)^3}\left((n^2+(1+2\alpha)n+\alpha)(1-x^2)-\alpha^2(1+x^2)+2\right),\quad 0<x<1,$$
then, it can be easily checked that for sufficiently large integer $n,$ we have $g'_{n,\alpha}(x)>0$ for $x\in [0, 1-\frac{\gamma}{n}].$ Hence, from Butlewski's theorem, see for example [\cite{Andrews}, p.238], the relative maxima of 
$|u_{\alpha}(x)|$ form a decreasing sequence. That is 
$$\sup_{x\in [0,1-\frac{\gamma}{n}]} (1-x^2)^{(1+\alpha)/2} |\qn(x)|\leq (1-t_{n,*}^2)^{(1+\alpha)/2} |\qn(t_{n,*})|,$$
where $t_{n,*}$ is the first zero  of $\big(Q_n^{(\alpha,\alpha)}\big)$ in $(0,1).$  Moreover, from the locations of the first two positive zeros of
$\qn(x)$ as well as from the local estimate of this later, given by \eqref{estimatesJacobi}, one concludes that  
$$(1-t_{n,*}^2)^{(1+\alpha)/2} |\qn(t_{n,*})|\leq C_{\alpha}.$$
On the other hand,  by using the previous inequality together with the local estimate of Jacobi polynomials, given by \eqref{estimatesJacobi}, one can easily check
that 
\begin{equation}
\label{Ineq2.3}
(1-y^2)^{1+\alpha}\left|\pn(x)\qn(y)- \pn(y)\qn(x)\right|\leq C_{\alpha} ,\quad\forall\,\, 0\leq x\leq y\leq 1-\frac{\gamma}{n}.
\end{equation}
Moreover, it has been shown in \cite{Area}, that 
\begin{equation}
\label{Ineq2.4}
\lim_{y\rightarrow 1} (1-y^2)^{\alpha} \qn(y) =\frac{1}{\sqrt{h_n^{(\alpha,\alpha)}}}\frac{2^{-1-\alpha}\Gamma(2+2\alpha)\Gamma(n+\alpha+1)}{\alpha\Gamma(1+\alpha)\Gamma(n+2\alpha+1)}\leq C_{\alpha} n^{-\alpha+1/2}.
\end{equation}
Also, it is well known that 
$$\sup_{x\in [0,1]}|\pn(x)| \leq C_{\alpha} n^{\alpha+\frac{1}{2}}.$$
Consequently, we have 
\begin{equation}
\label{Ineq2.5}
(1-y^2) |\pn(x)|\leq C_{\alpha} \left((1-(1-\frac{\gamma}{n})^2\right)n^{\alpha+\frac{1}{2}}\leq  C_{\alpha}n^{\alpha-\frac{1}{2}},\quad 1-\frac{\gamma}{n}\leq y\leq 1,\,\, x\in [0, y].
\end{equation}
Hence, by combining \eqref{Ineq2.3}, \eqref{Ineq2.4} and \eqref{Ineq2.5}, one gets 
\begin{equation}
\label{Ineq2.6}
(1-y^2)^{1+\alpha} |k_n(x,y|\leq C_{\alpha},\quad \forall\,\,\,  0\leq x\leq y\leq 1,
\end{equation}
where $k_n(x,y)$ is as given by \eqref{Eq3}. By using \eqref{Eq4}, \eqref{Wronskian} and the previous inequality, one 
gets the following bound for the reminder term $R_{n,\alpha}(x),$ given by \eqref{Eq2},
\begin{equation}
\label{Ineq2.7}
|R_{n,\alpha}(x)|\leq C_{\alpha} \frac{c^2}{2n+2\alpha+1} \int_{x}^1 |\ps(t)|\, dt.
\end{equation}
On the other hand, since $0<\alpha<3/2,$ then by using  the notations and the results of Theorem 1 and Proposition 2  and by using \eqref{boundJ}, one 
concludes that for sufficiently large values of $n,$ with $q=c^2/\chi_{n,\alpha}\leq q_0<1,$ we  have  for $0\leq x\leq 1,$ 
\begin{equation}
\label{Ineq2.8}
\int_x^1 |\ps(t)|\, dt \leq  \frac{C_{\alpha}(q_0)}{(1-q_0)^{1/4}} c_{\alpha} \int_0^1 (1-t)^{-\alpha/2-1/4}\, dt +\int_0^1
|\mathcal{E}_{n,\alpha}(t)| \, dt \leq C_{\alpha}(q_0).
\end{equation}
The previous two inequalities imply that 
\begin{equation}
\label{Ineq2.9}
\sup_{x\in [0,1]} |R_{n,\alpha}(x)| \leq C_{\alpha}(q_0) \frac{c^2}{2n+2\alpha+1}.
\end{equation}
Also, since
$\pn$ is bounded on $[0,1]$ which is not the case for $\qn,$ then we have $B=0.$ This implies that 
$$\ps(x)= A_n \pn(x)+ R_{n,\alpha}(x).$$
Recall that $\ps$ and $\pn$ are normalized so they have a unit $L^2(I,\omega_{\alpha})-$norms.  Hence,
by using  \eqref{Ineq2.9}, one gets
$$ |1-A_n|=\left| \|\ps\|_{L^2(I,\omega_{\alpha})}-A_n\|\pn\|_{L^2(I,\omega_{\alpha})}\right|\leq \| R_n\|_{L^2(I,\omega_{\alpha})} \leq C_{\alpha}(q_0) \frac{c^2}{2n+2\alpha+1}.$$
This concludes the proof of the theorem.$\qquad \Box$

\section{Decay rate of the eigenvalues of the weighted finite Fourier transform operator.}

In this section, we give a precise super-exponential decay rate of the eigenvalues 
$\lambda_n^{(\alpha)}$ of the  operator ${\displaystyle \mathcal Q_c^{\alpha}=\frac{c}{2\pi}
\mathcal F_c^{{\alpha}^*}  \mathcal F_c^{\alpha},}$ which implies the  decay rate 
of the eigenvalues $\mu_n^{(\alpha)}$ of   the operator $\mathcal F_c^{\alpha}.$ The study of this  decay rate 
is done under the condition that $0<\alpha < 3/2$ and it is based on the uniform asymptotic approximations of the GPSWFs, given by the previous two sections.
It has been shown in \cite{Karoui-Souabni1}, that for $\alpha \geq 0,$ the sequence of the eigenvalues 
$\lambda_n^{(\alpha)}(c),$ arranged in the decreasing order $1 > \lambda_0(c) > \lambda_0(c)> \cdots > \lambda_n(c) >\cdots >0,$ satisfies the following monotonicity property with respect to the parameter $\alpha,$
$ 
\lambda_n^{(\alpha)}(c) \leq \lambda_n^{(\alpha')}(c),\quad \forall\, \alpha \geq \alpha' \geq 0.$
Moreover, it has been shown in \cite{Bonami-Karoui1} that in the special case where $\alpha=0,$ the eigenvalues 
$\lambda_n(c)=\lambda_n^{(0)}(c)$ decay asymptotically faster than  $e^{-2n\log\left(\frac{a n}{c}\right)}$ for any positive real number $0<a < \frac{4}{e}.$ The constant  $\frac{4}{e}$ is optimal in the sense that it cannot be replaced by a larger constant. As a consequence of the previous monotonicity property, one concludes that for $\alpha >0,$ the eigenvalues $\lambda_n^{(\alpha)}(c),$ decay also  faster than $e^{-2n\log\left(\frac{a n}{c}\right)}, \, 0< a < \frac{4}{e}.$
Also, note that in \cite{Wang2}, the authors have given the following  explicit formula 
for the eigenvalues $\mu_n^{(\alpha)}(c),$
\begin{equation}\label{mu_n}
\mu_n^{(\alpha)}(c) = i^n \sqrt{\pi} \frac{\Gamma(n+\alpha+1)\Gamma(n+2\alpha+1)}{\Gamma(n+\alpha+3/2)\Gamma(2n+2\alpha+1)} c^n \exp(\Phi_n^{(\alpha)}(c)),\quad \Phi_n^{(\alpha)}(c)=\int_0^c \frac{F_n(\tau,\alpha)-n}{\tau}\, d\tau,\quad c>0,
\end{equation}
where
\begin{equation}\label{F_nc}
F_n(c,\alpha)= \int_{-1}^1 x \, \ps(x) \, \partial_x \ps(x)\, \omega_{\alpha}(x) \, dx.
\end{equation}
Hence, under the condition that the quantity $\Phi_n(c)$ is bounded and by using the bounds 
of the $\Gamma(\cdot),$ given by \eqref{Ineq2}, one gets the following super-exponential decay rate of 
the ${\displaystyle \lambda_n^{(\alpha)}(c)=\frac{c}{2\pi} |\mu_n^{(\alpha)}(c)|^2,}$
\begin{equation}\label{lambda_n}
\lambda_n^{(\alpha)}(c) \leq C_{\alpha} e^{-(2n+1) \log\left(\frac{4n+4\alpha+2}{e c}\right)} e^{\Phi_n^{(\alpha)}(c)},
\end{equation}
for some constant $C_{\alpha}$ and for large enough values of the integer $n.$ Note that comparing to our notations,
the roles of 
$\lambda_n^{(\alpha)}(c)$ and $\mu_n^{(\alpha)}(c)$ are reversed in  \cite{Wang2}. Also, in \cite{Wang2}, the authors have shown the convergence of the quantity $\Phi_n^{(\alpha)}(c)$ under the strong assumptions that 
$\ps(x)$ and $\partial_x \ps$ are well approximated by their projections over the five dimensional subspaces 
$\mbox{Span} \{\widetilde P_{n+2k}^{(\alpha,\alpha)}(x),\,\, -2\leq k\leq 2\},$ $\mbox{Span} \{\partial_x \widetilde P_{n+2k}^{(\alpha,\alpha)}(x),\,\, -2\leq k\leq 2\},$ respectively. Also, the given proof is based on the following equality, 
\begin{equation}\label{F_n0}
F_n(0,\alpha)= \int_{-1}^1 x \, \pn(x) \, \partial_x \pn(x)\, \omega_{\alpha}(x) \, dx = n.
\end{equation}
In the sequel, we prove the super-exponential decay rate of the 
$\lambda_n^{(\alpha)}(c)$ with $0<\alpha < \frac{3}{2}.$ This is given by the following proposition.

\begin{proposition} Let $c>0$ and $0<\alpha <\frac{3}{2}$ be two positive real numbers. Then, there exists 
$N_{\alpha}(c)$ and a constant $C_{\alpha}>0$ such that 
\begin{equation}
\label{decay_lambdan}
\lambda_n^{(\alpha)}(c) \leq C_{\alpha} \exp\left(-(2n+1)\left[ \log\left(\frac{4n+4\alpha+2}{e c}\right)+ C_{\alpha} \frac{c^2}{2n+1}\right]\right),\quad 
\forall\, n\geq N_{\alpha}(c).
\end{equation}
\end{proposition}

\noindent
{\bf Proof:} We recall that $C_{\alpha}$ is a generic constant that might take different values.
We choose $N_{\alpha}(c)\in \mathbf N,$ large enough so that the conditions of Theorem 2 are satisfied,
whenever $n\geq N_{\alpha}(c).$ Also, we let $C_{\alpha,q_0}= C_{\alpha}.$
Since $\alpha>0,$ and since $\omega_{\alpha}(\pm 1)=0,$ then by using an integration by parts, 
we rewrite the quantity $F_n(c,\alpha)$ as follows,
\begin{eqnarray*}
F_n(c,\alpha)&=& -\frac{1}{2} \int_{-1}^1 (\ps(x))^2\, \omega_{\alpha}(x)\, dx +\alpha \int_{-1}^1 (\ps(x))^2
x^2 (1-x^2)^{\alpha-1}\, dx\\
&=& -\frac{1}{2}-\alpha +\alpha \int_{-1}^1 (\ps(x))^2 \omega_{\alpha-1}(x)\, dx.
\end{eqnarray*}
Note that by replacing $\ps(x)$ by $\pn(x)$ in the previous equality and by using \eqref{F_n0}, one gets the identity
\begin{equation}\label{FF_n0}
\alpha \int_{-1}^1  (\pn(x))^2 \, \omega_{\alpha-1}(x) \, dx = n+\alpha+\frac{1}{2},\quad \alpha >0.
\end{equation}
Moreover, from Theorem 2, we have 
$$\ps(x)= A_n \pn(x)+R_{n,\alpha}(x),$$ where
\begin{equation}
\label{bounds}
|1- A_n|\leq C_{\alpha} \frac{c^2}{2n+2\alpha+1},\qquad \sup_{x\in [0,1]} |R_{n,\alpha}(x)|\leq C_{\alpha} \frac{c^2}{2n+2\alpha+1}.
\end{equation}
Hence, we have
\begin{eqnarray*}
\alpha \left|(\ps(x))^2-(\pn(x))^2\right|&\leq& |1-A_n^2| \alpha (\pn(x))^2\\
&&\quad +(|1-A_n|+|1+A_n|)  |R_{n,\alpha}(x)|  \alpha |\pn(x)|+ \alpha |R_{n,\alpha}(x)|^2.
\end{eqnarray*} 
Since from \eqref{bounds}, we have $|1-A_n^2|\leq C_{\alpha}\frac{c^2}{2n+2\alpha+1},$ then by using
\eqref{FF_n0}, one concludes that $${\displaystyle |1-A_n^2| \alpha \int_{-1}^1 (\pn(x))^2\, \omega_{\alpha-1}(x)\, dx\leq C_{\alpha} c^2.}$$
Also, by using  \eqref{FF_n0}, \eqref{bounds} and H\"older's inequality, one gets
\begin{eqnarray*}
\int_{-1}^1 |R_{n,\alpha}(x)|  \alpha |\pn(x)|\, \omega_{\alpha-1}(x)\, dx&\leq &C_{\alpha}
\frac{c^2}{2n+2\alpha+1} \alpha \left(\int_{-1}^1|\pn(x)|\, \omega_{\alpha-1}(x)\, dx\right)^{1/2}\\
&\leq & C_{\alpha}
\frac{c^2}{\sqrt{2n+2\alpha+1}}.
\end{eqnarray*} 
Similarly, by using  \eqref{bounds}, one gets ${\displaystyle \alpha \int_{-1}^1 |R_{n,\alpha}(x)|^2 \, \omega_{\alpha-1}(x)\, dx \leq C_{\alpha} \left(\frac{c^2}{2n+2\alpha+1}\right)^2.}$ By collecting everything together, one concludes that
\begin{eqnarray*}
\left|F_{n}(\alpha,\tau)-n\right|&=& \alpha \int_{-1}^1 \left((\psi_{n,\tau}^{(\alpha)}(x))^2-(\pn(x))^2 \right)\,\omega_{\alpha-1}(x)\, dx\leq C_{\alpha} \tau^2.
\end{eqnarray*}
Consequently, we have 
$$ \Phi_n^{(\alpha)}(c) = \int_0^c \frac{F_n(\tau,\alpha)-n}{\tau}\, d\tau \leq C_{\alpha} c^2.$$
To conclude for the proof of \eqref{decay_lambdan}, il suffices to combine the previous inequality and 
\eqref{lambda_n}.\\

Next, we give an asymptotic lower bound for the counting number of the eigenvalues $\lambda_n^{(\alpha)}(c).$
To this end, we first recall that 
\begin{equation}\label{Qc}
Q_c^{(\alpha)}(x)=\frac{c}{2\pi}\int\limits_{-1}^{1}K_{\alpha}(c(x-y))g(y)\omega_{\alpha}(y)dy
\end{equation}
Where 
\begin{equation} \label{noyau}
 K_{\alpha}(x,t)=\frac{c}{2\pi}\sqrt{\pi}2^{\alpha+\frac{1}{2}}\Gamma(\alpha+1)\frac{J_{\alpha+1/2}(c(x-t))}{(c(x-t))^{\alpha+1/2}}
\end{equation}
and the $\lambda_n^{(\alpha)} $ are the  eigenvalues of $Q_c^{\alpha}$  arranged in decreasing order.
\begin{theorem}
 Let $ 0<\delta<1$ and let $M_c(\delta)$ be the number of eigenvalues of $Q_c^{(\alpha)}$ , $ \alpha>0,$ which are not smaller than $ \delta.$ Then, we have
  \begin{equation}\label{IneqMc}
 \frac{\gamma_{\alpha}-\delta}{1-\delta}\frac{c}{2\pi} \left(2^{2\alpha+1}B(\alpha+1,\alpha+1)\right)^2 +o(c)\leq  M_c(\delta) \leq \frac{1}{\delta}\Bigg[ \frac{c}{2\pi} \Big[ 2^{2\alpha+1} B(\alpha+1,\alpha+1) \Big]^2 \Bigg]
  \end{equation}
 Where $ \gamma_{\alpha}=2^{4\alpha}\Big( \frac{B(2\alpha+1,2\alpha+1)}{B(\alpha+1,\alpha+1)} \Big)$ and $B(\cdot,\cdot)$ is the beta function
 \end{theorem}
 
\noindent
{\bf Proof:}  To obtain the lower bound estimate of $M_c(\delta),$ we use Marzo's formula, see \cite{Marzo} or \cite{Abreu}, 
 \begin{equation}\label{marzo}
  M_c(\delta)\geq Trace(Q_c^{(\alpha)}) - \frac{1}{1-\delta}(Trace(Q_c^{(\alpha)})-Norm(Q_c^{(\alpha)}))
 \end{equation}
 Note that the $ Trace(Q_c^{\alpha})$ has been already given in \cite{Wang2}, where it has been shown that 
\begin{equation}
\frac{2\pi}{c}\sum_{n=0}^{\infty}\lambda_n^{(\alpha)}(c)=\pi\frac{\Gamma^2(\alpha+1)}{\Gamma^2(\alpha+\frac{3}{2})}
\end{equation}
 Moreover, by using the well known identity  $
  \frac{\Gamma(\alpha+1)}{\Gamma(\alpha+3/2)}=\frac{(\Gamma(\alpha+1))^2 2^{2\alpha+1}}{\Gamma(2\alpha+2)\sqrt{\pi}}
  $ one gets 
 \begin{equation}\label{trace}
 Trace(Q_c^{\alpha})=\sum_n \lambda_n^{(\alpha)}=\frac{c}{2\pi} \Big[ 2^{2\alpha+1} B(\alpha+1,\alpha+1) \Big]^2.
 \end{equation}
 To compute an estimate of $ Norm(Q_c^{(\alpha)}),$ we proceed as follows :
  $$Norm(Q_c^{\alpha}) =\sum_{n=0}^{\infty}(\lambda_n^{(\alpha)})^2=\frac{c^2}{4\pi} 2^{2\alpha +1}\Gamma^2(\alpha+1) \int\limits_{-1}^{1} \int\limits_{-1}^{1} \frac{J_{\alpha+1/2}^2(c(x-y))}{(c(x-y))^{2\alpha+1}} (1-x^2)^{\alpha}(1-y^2)^{\alpha}dxdy $$
We apply the change of variable $ y=\sigma $ and $ x=\sigma+\frac{\tau}{c} $, to obtain :
\begin{equation}\label{norm1}
 \sum_{n=0}^{\infty}(\lambda_n^{(\alpha)})^2=\frac{c}{4\pi} 2^{2\alpha +1}\Gamma^2(\alpha+1)\int_{-1}^{1}(1-\sigma^2)^{\alpha} \int_{c(-1+\sigma)}^{c(1+\sigma)} \frac{J^2_{\alpha+1/2}(\tau)}{\tau^{2\alpha+1}}\Big(1-(\sigma+\frac{\tau}{c})^2 \Big)^{\alpha} d\tau d\sigma.
\end{equation}
Since $$ \frac{J^2_{\alpha+1/2}(\tau)}{\tau^{2\alpha+1}}(1-\sigma^2)^{\alpha}\Big(1-(\sigma+\frac{\tau}{c})^2 \Big)^{\alpha} \leq \frac{J^2_{\alpha+1/2}(\tau)}{\tau^{2\alpha+1}} (1-\sigma^2)^{\alpha} \left(1-(\sigma-\frac{1}{c})^2\right)^{\alpha} $$
and since from [\cite{NIST},  p. 244], we have
\begin{equation}\label{eq1}
\int_{\mathbb{R}} \frac{J^2_{\alpha+1/2}(\tau)}{\tau^{2\alpha+1}} d\tau = \frac{1}{2^{2\alpha}}\Gamma(1/2)\frac{\Gamma(2\alpha+1)}{\Gamma(2\alpha+3/2)\Gamma^2(\alpha+1)}
\end{equation}
 then by Lebesgue's dominated convergence theorem applied to the  integral in (\ref{norm1}), one gets 
 \begin{equation*}
 \lim_{c\to \infty }\Big[ \frac{1}{c} \sum_n(\lambda_n^{(\alpha)})^2 \Big] = \frac{1}{4\pi}2^{2\alpha+1}\Gamma^2(\alpha+1) \int_{-1}^{1}(1-\sigma^2)^{2\alpha}d\sigma \int_{\mathbb{R}}\frac{J^2_{\alpha+1/2}(\tau)}{\tau^{2\alpha+1}} d\tau.
 \end{equation*}
 Since $ \int_{-1}^{1}(1-\sigma^2)^{2\alpha}d\sigma = 2^{4\alpha+1}B(2\alpha+1,2\alpha+1) $ , then by straightforward computations, one gets 
 $$
 \lim_{c\to \infty }\Big[ \frac{1}{c} \sum_n(\lambda_n^{(\alpha)})^2 \Big] = \gamma_{\alpha}\frac{1}{2\pi} \Big[ 2^{2\alpha+1} B(\alpha+1,\alpha+1) \Big]^2,\quad \gamma_{\alpha}= 2^{4\alpha} \Big( \frac{B(2\alpha+1,2\alpha+1)}{B(\alpha+1,\alpha+1)} \Big)^2.
 $$
Here, $B(\cdot,\cdot)$ is the Beta function.  Hence, for $\alpha>-1,$ we have  
   \begin{equation} \label{norme}
   \sum_n(\lambda_n^{(\alpha)})^2 =\gamma_{\alpha}\frac{c}{2\pi} \Big[ 2^{2\alpha+1} B(\alpha+1,\alpha+1) \Big]^2+o(c),
  \end{equation}
  To conclude for the proof of the lower bound estimate in \eqref{IneqMc}, it suffices to combine  (\ref{norme}) and (\ref{trace}) in (\ref{marzo}). 
  
   Finally, to prove the upper bound in \eqref{IneqMc},
  it suffices to note that 
  $$
   \sum_{k=0}^{\infty}\lambda_k^{(\alpha)} \geq \sum_{k=0}^{M_c(\delta)}\lambda_k^{(\alpha)} \geq \delta M_c(\delta) 
$$ and then use \eqref{trace}.$\qquad \Box$

 \begin{remark}
  In the special case $ \alpha=0,$ the inequalities \eqref{IneqMc} become  
   \begin{equation}
   \frac{2}{\pi}+o(1) \leq \frac{M_c(\delta)}{c}\leq \frac{2}{\pi \delta}
   \end{equation}
In the special case $\alpha=0,$ it  has been shown in \cite{Landau2} that $ M_c(\delta) $ is independent of $ \delta\in [0,1].$ Hence, for $\alpha=0$ and by  letting $\delta \rightarrow 1$ in \eqref{IneqMc},
 we recover the following  Landau's classical result, see \cite{Landau2}, 
   $M_c(\delta)=\frac{2c}{\pi}+o(c).$
    \end{remark}

\end{document}